\documentclass[final, 3p, fleqn, authoryear]{elsarticle}
\usepackage{graphicx}[dvipdfmx]
\makeatletter
\renewcommand{\S}[1]{#1^{\mathrm{S}}}
\makeatother

\makeatletter
\def\ps@pprintTitle{%
 \let\@oddhead\@empty
 \let\@evenhead\@empty
 \def\@oddfoot{}%
 \let\@evenfoot\@oddfoot}
\makeatother

\usepackage{chngcntr}
\usepackage{apptools}
\AtAppendix{\counterwithin{lem}{section}}
\usepackage{amsmath,amssymb}
\usepackage{mathtools}
\usepackage[ISO]{diffcoeff}
\usepackage{rmss_symbol}
\usepackage{amsthm}
\theoremstyle{definition}
\newtheorem{thm}{Theorem}[section]

\newtheorem{cor}{Corollary}

\newtheorem{lem}{Lemma}[section]

\usepackage{refcount}
\ExplSyntaxOn
\NewDocumentCommand{\replace}{mmm}
 {
  \marian_replace:nnn {#1} {#2} {#3}
 }

\tl_new:N \l_marian_input_text_tl
\tl_new:N \l_marian_search_tl
\tl_new:N \l_marian_replace_tl

\cs_new_protected:Npn \marian_replace:nnn #1 #2 #3
 {
  \tl_set:Nn \l_marian_input_text_tl { #1 }
  \tl_set:Nn \l_marian_search_tl { #2 }
  \tl_set:Nn \l_marian_replace_tl { #3 }
  \regex_replace_all:nnN { \b\u{l_marian_search_tl}\b } { \u{l_marian_replace_tl} } \l_marian_input_text_tl
  \tl_use:N \l_marian_input_text_tl
 }
\ExplSyntaxOff

\newcommand{\thmref}[1]{Non-existence theorem}

\usepackage{chngcntr}
\usepackage{apptools}
\AtAppendix{\counterwithin{lem}{section}}

\usepackage{nameref}

\usepackage{txfonts}
\usepackage{pxfonts}

\usepackage{color}

\usepackage{float}
\restylefloat{table}

\usepackage{hyperref}
\hypersetup{pdfauthor={Name},colorlinks=true,urlcolor=blue}

\journal{{}}

\biboptions{authoryear}

\begin{document}

\begin{frontmatter}

\title{
Non-existence of queues
\\
for system optimal departure patterns in tree networks}

\author[1]{Takara Sakai\corref{cor1}}
\ead{takara.sakai.t1@dc.tohoku.ac.jp}
\author[1]{Koki Satsukawa\corref{cor1}}
\ead{satsukawa@tohoku.ac.jp}
\author[1]{Takashi Akamatsu}
\ead{akamatsu@plan.civil.tohoku.ac.jp}

\address[1]{Graduate School of Information Sciences, Tohoku University, 6-6 Aramaki Aoba, Aoba-ku,
Sendai, Miyagi 980-8579, Japan}

\cortext[cor1]{Corresponding author}

\begin{abstract}
This study proves the non-existence of queues for a dynamic system optimal (DSO) departure pattern in a directed rooted tree network with a single destination.
First, considering queueing conditions explicitly, we formulate the DSO problem as mathematical programming with complementarity constraints (MPCC) that minimizes the total system cost which consists of the schedule and queueing delay costs.
Next, for an arbitrary feasible solution to the MPCC, we prove the existence of another feasible solution where the departure flow pattern on every link is the same but no queue exists.
This means that the queues can be eliminated without changing the total schedule delay cost. 
Queues are deadweight losses, and thus the non-existence theorem of queues in the DSO solution is established.
Moreover, as an application of the non-existence theorem, we show that the MPCC can be transformed into a linear programming (LP) problem by eliminating the queueing conditions.
\end{abstract}

\begin{keyword}
	\textit{
		departure time pattern, rooted tree network, bottleneck model, 
		dynamic system optimal assignment, global optimal solution
	}
\end{keyword}

\end{frontmatter}

\section{Introduction}
In Vickrey's bottleneck model~\citep{Vickrey1969-ic}, a system optimal departure pattern, which minimizes the total system cost, does not cause a bottleneck queue. 
This fact makes intuitive sense in the case of a single bottleneck. 
However, it is not so apparent in more general settings such as networks with multiple bottlenecks or with heterogeneous commuters.
In other words, we cannot immediately answer the question of whether eliminating queues always improves the system's efficiency.
In fact, to the best of our knowledge, there are no studies that answer this question.

\par
In this paper, as a preliminary step to answering the above question for more general settings, we investigate the dynamic system optimal (DSO) problem in a directed rooted tree network with a many-to-one origin-destination pattern. 
We first formulate the DSO problem straightforwardly as mathematical programming with complementarity constraints (MPCC) that minimizes the sum of the schedule and queueing delay costs (total system cost) with queueing conditions.
We then show that for an arbitrary departure pattern with queues, we can construct another departure pattern where no queues exist without changing the total schedule delay costs.
This proves that the queueing delay costs are deadweight loss, and thus they must not exist in a DSO state.
Moreover, as a mathematical application of this fact, we show that the DSO problem can be formulated as a linear programming (LP) problem.
In other words, we find a possibility of a transformation from an MPCC to an LP under a certain condition.

\par
The remainder of this paper is organized as follows. 
Section 2 describes the model's settings. 
Section 3 formulates the DSO problem. 
Then, in Section 4, we prove the non-existence theorem and show some interesting results derived from the theorem. 
Section 5 concludes the paper.

\begin{figure}[tbp]
		\center
		\includegraphics[clip, width=0.6\columnwidth]{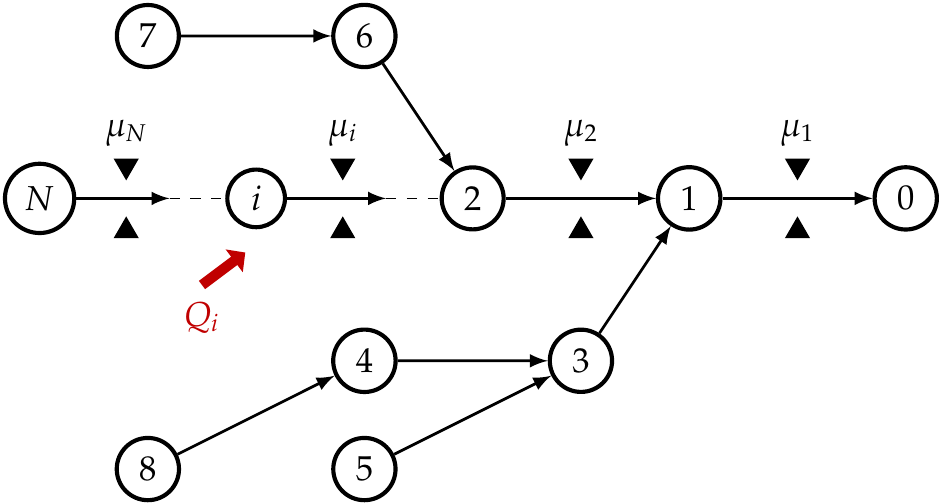}
		\caption{Directed rooted tree network with a many-to-one origin destination pattern; every link has a bottleneck}
    \label{Fig_TN}
\end{figure}

\section{Model settings}
Consider a directed rooted tree network consisting of $N$ origin nodes and a single destination node, as shown in Figure \ref{Fig_TN}.
The destination is the root and numbered $0$.
The origins are appropriately numbered from $1$ to $N$, and the set of origins is denoted by $\mathcal{N}\equiv \{i \mid i=1,2,.... N\}$.
We denote by $\underline{\mathcal{N}}(i)$ the set of an origin $i$ and the (downstream) ascendants of the origin $i$; we also denote by $\overline{\mathcal{N}}(i)$ the set of an origin $i$ and the (upstream) descendants of the origin $i$.

The edge (link) connecting from an origin $i\in\mathcal{N}$ to its parent is referred to as link $i$.
Hence, we utilize the notation $\mathcal{N}$ interdependently, this indicates the set of nodes or links depending on the context.
A link $i\in\mathcal{N}$ has a bottleneck with a finite capacity $\mu_{i}$ at the end of the link.
For the pattern of the bottleneck capacity, we assume that the following condition is satisfied:
\begin{align}
    &
    \sum_{j\in\mathcal{B}(i)} \mu_{j}\leq \mu_{i},
    && \forall i\in\mathcal{N},
		\label{Eq_Asm_Capa}
\end{align}
where $\mathcal{B}(i)$ is the set of incoming links to node $i$. 

At each bottleneck, a queue is formed when the inflow flow rate exceeds the capacity.
The queue evolution and associated queueing delay are modeled by the standard point queue model in accordance with the first-in-first-out (FIFO) principle.
Here, we denote by $A_{i}(t)$ and $D_{i}(t)$ the cumulative arrival and departure flows by time $t$ at a link $i$, respectively.
These cumulative flows, which can be regarded as functions of time $t$, are  assumed to be Lipschitz continuous; this guarantees that the arrival flow and departure flow rates are finite values and physically plausible.
The free flow travel time of a link $i$ is denoted by $d_{i}$.
We also denote by $\underline{i}$ the nearest downstream link of link $i$.

\begin{figure}
	\center
	\includegraphics[clip, width=0.55\columnwidth]{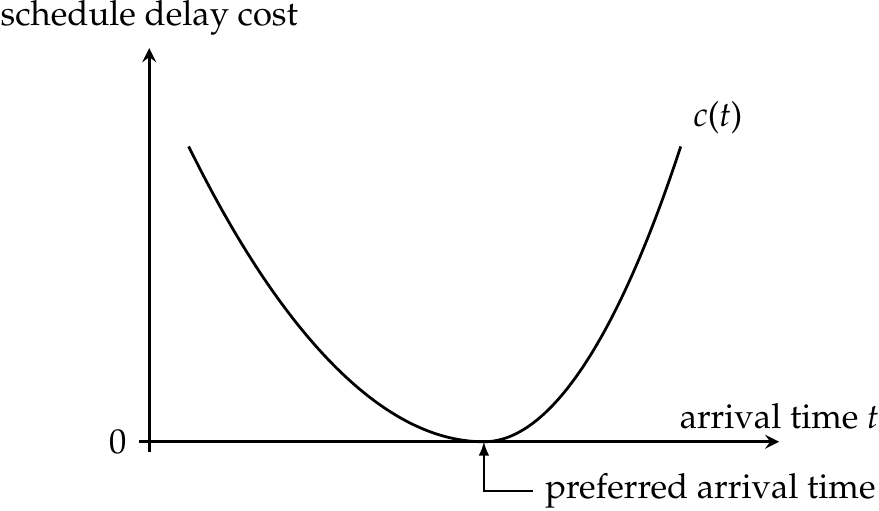}
	\caption{Schedule delay cost functions}
	\label{Fig_SDF}
\end{figure}

The cost of the commuters' trip is assumed to be separable into free-flow travel, queueing delay and schedule delay costs.
A schedule delay cost is associated with the difference between the actual and preferred arrival times at the destination, as shown in Figure \ref{Fig_SDF}.
We assume that commuters have the same preferred arrival time and the same value for a schedule delay.
They have the same schedule delay cost function $c(t): \ClT \rightarrow \mathbb{R}_{+}$ for an actual destination arrival time $t$.
We also assume that the value of $c(t)$ is zero if $t$ is the preferred arrival time for each group.

\section{Dynamic system optimal assignment problem}

\subsection{Formulation in a Lagrangian-like coordinate system}
In modeling a traffic assignment problem, we describe traffic variables mainly in a \textit{Lagrangian-like coordinate system}~\citep[][]{Akamatsu2015-ip,Akamatsu2021-zg}.
In this system, variables are expressed as functions of the arrival time at the destination, and not at the origin or bottleneck. 
Such an expression is suitable for considering the ex-post travel time of each commuter during his/her trip. 
Hence, we can easily trace the time-space paths of commuters, which allows for the representation of the commuters' trip costs in a simple manner.

\begin{figure}[tbp]
	\center
	\includegraphics[clip, width=0.80\columnwidth]{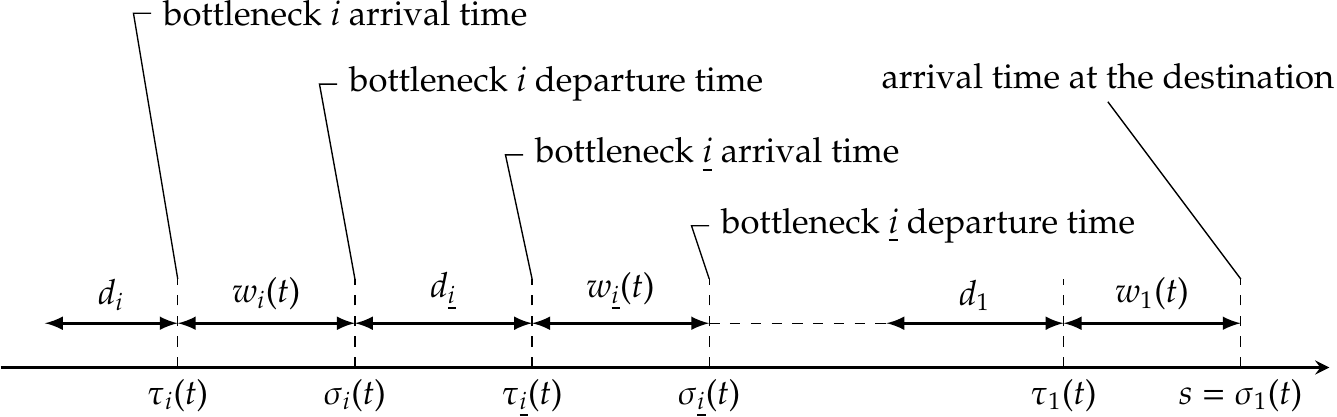}
	\caption{Time-space path of commuters who arrive at the destination at time $t$}
  \label{Fig_TS-path}
\end{figure}

Based on the concept of the Lagrangian-like coordinate system, we define the arrival times $\tau_{i}(t)$ and departure times $\sigma_{i}(t)$ at a bottleneck $i\in\mathcal{N}$ for commuters whose destination arrival time is $t$.
Mathematically, these variables are represented as follows (see also Figure~\ref{Fig_TS-path}):
\begin{align}
	&\tau_{i}(t) = t 
	- \sum_{j \in \underline{\ClN}(i)} w_{j}(t) 
	- \sum_{j \in \underline{\ClN}(i) \setminus \{ i \}} d_{j}
	&&\forall i \in \ClN, \quad \forall t \in \ClT,
	\label{Eq_tau=s-sumw}
	\\
	&\sigma_{i}(t) = t 
	- \sum_{j \in \underline{\ClN}(i) \setminus \{ i \}} w_{j}(t)
	- \sum_{j \in \underline{\ClN}(i) \setminus \{ i \}} d_{j} 
	&&\forall i \in \ClN, \quad \forall t \in \ClT,
	\label{Eq_sigma=s-sumw}
\end{align}
where $w_{i}(t)$ is the queuing delay at a bottleneck $i$ for the commuters whose destination arrival time is $t$.
Note that the derivative satisfies the following relationship:
\begin{align}
    &\cfrac{\mathrm{d}\tau_{i}(t)}{\mathrm{d}t}\equiv \dot{\tau}_{i}(t) > 0,
    &&\forall i \in \ClN, \quad \forall t \in \ClT,
    \label{Eq:TauMonotone}
\end{align}
where an overdot of a variable means the derivative of the variable.
This is derived from the assumption of the Lipschitz continuity of the cumulative flows\footnote{
The Lipschitz continuity requires that the following relationship holds:
\begin{align*}
 &0 \leq \diff{ A_{i}(\tau_{i}(t))}{\tau_{i}(t)} < \infty
&&\forall i \in \ClN, \quad \forall t \in \ClT.
\label{Eq_Lipschitz}
\end{align*}
Then, the chain rule derives $\mathrm{d}A_{i}(\tau_{i}(t))/\mathrm{d}\tau_{i}(t) = (\mathrm{d}A_{i}(\tau_{i}(t))/\mathrm{d}t) (1/\dot{\tau}_{i}(t))$.
It is obvious that $\mathrm{d}A_{i}( \tau_{i}(t))/\mathrm{d}t\geq 0$; otherwise, the cumulative curve becomes \textit{backward-bending}, which indicates the existence of a negative flow. Hence, $\dot{\tau}_{i}(t) > 0$ must hold.}.

\subsection{Physical conditions in a traffic assignment problem}
We formulate the physical conditions for a traffic assignment problem.
First, the inflow from each origin must satisfy the following demand conservation condition:
\begin{align}
    &\int_{\ClT} q_{i}(t) \mathrm{d} t = Q_{i}
	&&\forall i \in \ClN,
\end{align}
where $q_{i}(t)\geq 0$ is the inflow rate to the network of $(i)$-commuters whose destination arrival time is $t$.

Next, we describe the queueing congestion at each bottleneck by using a point queue model, as mentioned in the previous section.
Let $x_{i}(t)$ be the departure flow rate at link $i\in\mathcal{N}$.
Then, the queueing condition is described as the following complementarity condition:
\begin{align}
    &\begin{dcases}
    x_{i}(\sigma_{i}(t))
    =
    \mu_{i}
    \quad &\mathrm{if}\quad
    w_{i}(t) > 0
    \\
    x_{i}(\sigma_{i}(t))
    \leq
    \mu_{i}
    \quad &\mathrm{if}\quad
    w_{i}(t) = 0
  \end{dcases}
	&&\forall i \in \ClN,
  \quad \forall t \in \ClT.
	\label{Eq_queueing_x}
\end{align}

Third, we formulate the flow conservation condition at each node.
This condition requires that for each time, the inflow rate to a node should equal the outflow rate from the node.
Mathematically, this is described as follows:
\begin{align}
    &\sum_{j\in\mathcal{B}(i)}\cfrac{\mathrm{d}D_{j}(\sigma_{j}(t))}{\mathrm{d}t} + q_{i}(t) = \cfrac{\mathrm{d}A_{i}(\tau_{i}(t))}{\mathrm{d}t}
    &&
    \quad \forall i\in\mathcal{N},
    \quad \forall t\in\mathcal{T},
\end{align}
where $\mathcal{B}(i)$ is the set of children of node $i$, of which $i$ is the parent.
Here, the FIFO principle is written as follows:
\begin{align}
    &A_{i}(\tau_{i}(t)) = D_{i}(\sigma_{i}(t)),
    &&
    \quad \forall i\in\mathcal{N},
    \quad \forall t\in\mathcal{T}.
\end{align}
Substituting this condition into the flow conservation condition, we derive the following condition in summary:
\begin{align}
    &\sum_{j\in\mathcal{B}(i)}x_{j}(\sigma_{j}(t))\dot{\sigma}_{j}(t) + q_{i}(t) =
    x_{i}(\sigma_{i}(t))\dot{\sigma}_{i}(t),
    &&
    \quad \forall i\in\mathcal{N},
    \quad \forall t\in\mathcal{T}.
\end{align}

It is worth noting that recursively applying this condition to the descendants of each node $i$, we have the following equivalent condition for the flow conservation condition:
\begin{align}
    \sum_{j\in\overline{\mathcal{N}}(i)} q_{i}(t) = x_{i}(\sigma_{i}(t))\dot{\sigma}_{i}(t),
    &&
    \quad \forall i\in\mathcal{N},
    \quad \forall t\in\mathcal{T}.\label{Eq:OutflowRelation}
\end{align}
Thus, substituting this into Eq.~\eqref{Eq_queueing_x}, we can describe the queueing condition by utilizing the variable $q_{i}(t)$, as follows:
\begin{align}
    &\begin{dcases}
    \sum_{j \in \overline{\ClN}(i)} q_{j}(t)
    = \mu_{i} \dot{\sigma}_{i}(t)
    \quad &\mathrm{if}\quad
    w_{i}(t) > 0
    \\
		\sum_{j \in \overline{\ClN}(i)} q_{j}(t)
    \leq \mu_{i} \dot{\sigma}_{i}(t)
    \quad &\mathrm{if}\quad
    w_{i}(t) = 0
   \end{dcases}
	 &&\forall i \in \ClN,
   \quad \forall t \in \ClT.
   \label{Eq:QueueCondRevised}
\end{align}

\subsection{Formulation of a dynamic system optimal assignment problem}
Under the abovementioned physical conditions, we formulate a dynamic system optimal (DSO) assignment problem that derives an optimal state minimizing the total system cost.
Here, the total system cost is defined as the sum of the commuters' trip costs, which consist of the queueing delay and schedule delay costs.
Mathematically, the trip cost experienced by $(i)$-commuters whose destination arrival time is $t$ is described as follows:
\begin{align}
  C_{i}(t) =
  c(t) +
  \alpha \sum_{j \in \underline{\ClN}(i)} \left( w_{j}(t) + d_{j} \right)
  &&\forall i \in \ClN,
  \label{Eq_TripCost}
\end{align}
where $\alpha$ is a parameter representing the value of time.
This study assumes $\alpha=1$ regardless of a commuter's group (i.e., the cost sensitivity to queueing delay time is the same for all commuters).

The DSO assignment problem is then formulated as follows:
\begin{align}
	\mbox{[DSO]}
	\qquad
    \min_{\Vtq, \Vtw}.
    \quad
    &
		f(\Vtq, \Vtw) = 
    \sum_{i \in \ClN}
    \int_{\ClT} 
    \left(
      c(t) + \sum_{j;j \in \underline{\ClN}(i)}  \left( w_{j}(t) + d_{j} \right) 
      \right)
    q_{i}(t) \mathrm{d} t
    \label{Eq:ObjectFunc}
    \\
    \mbox{s.t.} \quad
    &\begin{dcases}
			\sum_{j \in \overline{\ClN}(i)} q_{j}(t)
			= \mu_{i} \dot{\sigma}_{i}(t)
			\quad &\mathrm{if}\quad
			w_{i}(t) > 0
			\\
			\sum_{j \in \overline{\ClN}(i)} q_{j}(t)
			\leq \mu_{i} \dot{\sigma}_{i}(t)
			\quad &\mathrm{if}\quad
			w_{i}(t) = 0
		 \end{dcases}
		 &&\forall i \in \ClN,
		 \quad \forall t \in \ClT.
		\label{Eq_DSO-MPCC_cnst_Queueing}
  \\
  &\int_{\ClT} q_{i}(t) = Q_{i}
  &&\forall i \in \ClN,
	\label{Eq_DSO-MPCC_cnst_Demand}
  \\
  &q_{i}(t) \geq 0
	&&\forall i \in \ClN,
  \quad \forall t \in \ClT,
	\label{Eq_DSO-MPCC_cnst_q_Nonneg}
	\\
  &\sum_{j \in \underline{\ClN}(i)} \dot{w}_{i}(t) < 1
	&&
	\forall i \in \ClN,
	\forall t \in \ClT,
	\label{Eq_DSO-MPCC_const_Tau}
	\\
 \mbox{where} \quad 
  &\dot{\sigma}_{i}(t) = 1 - \sum_{j \in \underline{\ClN}(i)} \dot{w}_{i}(t)
  &&\forall i \in \ClN, \quad \forall t \in \ClT.
\end{align}
The objective function is the total system cost, which is the sum of the trip costs \eqref{Eq_TripCost} of all the commuters.
The constraint \eqref{Eq_DSO-MPCC_const_Tau} is derived from Eq.~\eqref{Eq:TauMonotone}.
We denote by $\{ \S{q}_{i}(t) \}$ and $\{ \S{w}_{i}(t) \}$ the optimal solution to [DSO].
Note that $\sigma_{i}(t)$ is a dependent variable (see Eq.~\eqref{Eq_sigma=s-sumw}), and thus, it is not necessary to treat this as a decision variable (i.e. solution).
We denote by $\mathcal{X}$ the set of feasible solutions $(\boldsymbol{q},\boldsymbol{w})$ to [DSO].

Since the constraints of [DSO] include complementarity conditions, this problem is classified as \textit{mathematical programming with complementarity constraints} (MPCCs). 
In general, it is difficult to clarify the analytical properties (e.g. existence and uniqueness) of this kind of problems~\citep[][]{Luo1996-in} since they do not have good mathematical properties, such as convexity.
It is also known that the optimality conditions are in a complex form, which means that it is difficult to understand the characteristics of the solution.

\section{There are no queues for the global optimal solution in the tree network}

\subsection{Non-existence theorem}
We show our main theorem that states the non-existence of queues in the DSO state, as follows:
\begin{thm}[\textit{Non-existence theorem}]
    There \textit{do not exist} queues for the system optimal solution to [DSO]:
    \begin{align}
      &\S{w}_{i}(t)=0
      &&\forall i \in \mathcal{N},
      \quad \forall t \in \ClT.
      \label{Eq_No-queue_DSO}
    \end{align}
    \label{Thm_No-queue_DSO}
\end{thm}
\noindent This means that the congestion externality should be perfectly eliminated in the DSO state.

The non-existence theorem is more than a simple contribution toward clarifying the characteristics of the system optimal state; from a mathematical point of view, it also establishes the equivalence between the MPCCs and a \textit{linear programming} (LP), as follows:
\begin{cor}[\textit{Equivalent linear programming}]
	The global optimal solution $\{\S{q}_{i}(t)\}$ to [DSO] is equivalent to the solution of the following LP:
	\begin{align}
		\mbox{[DSO-LP]} \qquad 
		\min_{\Vtq}.
		\quad
		&
		\sum_{i \in \ClN}
		\int_{\ClT} 
			\left( c(t) + \sum_{j \in \underline{\ClN}(i)} d_{j} \right) q_{i}(t) \mathrm{d} t
		\\
		\mbox{s.t.} \quad
		&\sum_{j \in \overline{\ClN}(i)} q_{j}(t) \leq \mu_{i} 
		&&\forall i \in \ClN, \quad \forall t \in \ClT,
		\label{Eq_DSO-LP_LinkFlow_BNCapa}
		\\
	&\int_{\ClT} q_{i}(t) = Q_{i}
	&&\forall i \in \ClN,
	\label{Eq_DSO-LP_ODcncv}
	\\
	&q_{i}(t) \geq 0
	&&\forall i \in \ClN,
	\quad \forall t \in \ClT.
	\label{Eq_DSO-LP_cnst_q_nonnega}
	\end{align}
\end{cor}
\noindent This corollary is straightforwardly derived from Eq.~\eqref{Eq_No-queue_DSO}.
We first see that the objective function becomes linear.
We next see that Constraint~\eqref{Eq_DSO-MPCC_cnst_Queueing} becomes the capacity constraint~\eqref{Eq_DSO-LP_LinkFlow_BNCapa}, which is linear inequality constraints.
Since queues do not exist, we can eliminate Constraint \eqref{Eq_DSO-MPCC_const_Tau}.
Thus, [DSO] reduces to [DSO-LP], which is a problem that minimizes the total schedule delay cost under the bottleneck capacity constraints.

This equivalence allows us to derive the following theoretical properties of the global optimal solution to the MPCCs from the standard mathematical optimization theory:
\begin{cor}[\textit{Existence}]
	Suppose that the assignment period $\ClT$ is sufficiently long and that all the bottleneck capacities are not zero. 
	Then, the solution to [DSO] always exists.
\end{cor}

\begin{cor}[\textit{Optimality condition}~\citep{Luenberger1997-la}]
	The optimal solution to [DSO] 
	($\{ \S{q}_{i}(t) \}$, $\{ \S{\rho}_{i} \}$, and 
	$\{ \S{p}_{i}(t) \}$) satisfies 
	the following optimality conditions: 
	\begin{align}
		\mbox{[DSO-LP-OC]} \qquad
		&  \int_{\ClT} \S{q}_{i}(t) \mathrm{d} t   = Q_{i}
		&&
		\forall i \in \mathcal{N},
		\\
		&\begin{dcases}
			\sum_{j \in \underline{\ClN}(i)} \left( \S{p}_{j}(t) + d_{j} \right) + c(t)  =  \S{\rho}_{i}
			\quad &\mathrm{if}\quad \S{q}_{i}(t) > 0
			\\
			\sum_{j \in \underline{\ClN}(i)}  \left( \S{p}_{j}(t) + d_{j} \right) + c(t)  \geq \S{\rho}_{i}
			\quad &\mathrm{if}\quad \S{q}_{i}(t) = 0
		\end{dcases}
		&&
		\forall i \in \mathcal{N},
		\quad \forall t \in \ClT,
		\label{Eq_DSO-OC_BNC}
		\\
		&\begin{dcases}
			\sum_{j \in \overline{\ClN}(i) } \S{q}_{j}(t) = \mu_{i}
			\quad &\mathrm{if}\quad \S{p}_{i}(t) > 0
			\\
			\sum_{j \in \overline{\ClN}(i) } \S{q}_{j}(t) \leq \mu_{i}
			\quad &\mathrm{if}\quad \S{p}_{i}(t) = 0
		\end{dcases}
		&&\forall i \in \mathcal{N},
		 \quad \forall t \in \ClT.
		\label{Eq_DSO-OC_DTC}
	\end{align}
	where $\rho_{i}$ and $p_{i}(t)$
	are Lagrangian multipliers for 
	the demand conservation condition \eqref{Eq_DSO-LP_ODcncv} 
	and the bottleneck capacity constraint \eqref{Eq_DSO-LP_LinkFlow_BNCapa},
	respectively. 
\end{cor}

\subsection{Proof}
We prove the non-existence theorem in a constructive manner.
Specifically, for an arbitrary traffic state, we devise a method to construct another traffic state where no queues exist and the pattern of link departure flow is the same, as follows:
\begin{lem}\label{Lemm:DeadWeight}
Consider a feasible traffic state $(\boldsymbol{q},\boldsymbol{w})\in\mathcal{X}$, and the pattern of link departure flow $\boldsymbol{x}$ resulting from Eq.~\eqref{Eq:OutflowRelation}.
Then, there exists another feasible traffic state $(\boldsymbol{q}^{*},\boldsymbol{w}^{*})\in\mathcal{X}$ satisfying the following relationship:
\begin{align}
    &w^{\ast}_{i}(t) = 0
	&&\forall i \in \ClN, \quad \forall t \in \ClT,
	\label{Eq:Trans_w}
	\\
	&x^{\ast}_{i}(\sigma_{i}^{*}(t)) = x_{i}(\sigma_{i}^{*}(t))
	&&\forall i \in \ClN, \quad \forall t \in \ClT,
	\label{Eq:Trans_x}
	\\
	&q^{\ast}_{i}(t) =
	x^{\ast}_{i}(\sigma^{\ast}_{i}(t)) - \sum_{j\in\mathcal{B}(i)}x^{\ast}_{j}(\sigma^{\ast}_{j}(t)),
	&&\forall i \in \ClN, \quad \forall t \in \ClT.
	\label{Eq:Trans_q}
\end{align}
\end{lem}
\noindent This lemma gives a sufficient condition for establishing the non-existence theorem.
Because the patterns of link departure flow $\boldsymbol{x}$ and $\boldsymbol{x}^{*}$ are the same, the sum of the schedule delay costs in the two traffic states does not change.
It follows that the total system cost in the new traffic state is always lower than that in the original traffic state.
In other words, this lemma claims that queueing delay costs are pure \textit{dead weight losses} in the tree network with homogeneous commuters.

The remainder of this section is mainly devoted to the proof of the lemma.
For this, it is sufficient to check whether the new traffic state is in a feasible region, i.e., it satisfies the conditions. 
Condition~\eqref{Eq_DSO-MPCC_const_Tau} is obviously satisfied since for an arbitrary node $i\in\mathcal{N}$, $w_{i}(t) = 0$ and thus $\dot{w}_{i}(t)=0<1$.
Therefore, we confirm that the other conditions are satisfied.

We first confirm that Condition~\eqref{Eq_DSO-MPCC_cnst_Queueing} holds.
Recursively applying Eq.~\eqref{Eq:Trans_q}, we have the following equation in the same manner as Eq.~\eqref{Eq:OutflowRelation}:
\begin{align}
    &
    x_{i}(\sigma_{i}^{*}(t)) = \sum_{j\in\overline{\mathcal{N}}(i)}q_{j}^{*}(t),
    &&
    \forall i\in\mathcal{N},\forall t\in\mathcal{T}.
\end{align}
We then have
\begin{align}
	&\sum_{j \in \overline{\ClN}(i)} q^{\ast}_{j}(t) - \mu_{i} \sigma^{\ast}_{i}(t)
	= x_{i}(\sigma^{\ast}_{i}(t)) - \mu_{i} \leq 0
	&&\because \sigma^{\ast}_{i}(t) = 1 \qquad 
	\forall i \in \ClN, \quad \forall t \in \ClT.
\end{align}
This is consistent with Condition~\eqref{Eq_DSO-MPCC_cnst_Queueing}.

Next, for an arbitrary node $i\in\mathcal{N}$, we can confirm that Condition~\eqref{Eq_DSO-MPCC_cnst_Demand} holds from the following simple calculation:
\begin{align}
    \int_{\ClT}
	q^{\ast}_{i}(t)
	\mathrm{d}t =
	\int_{\ClT} 
	\left(
		x_{i}(\sigma^{\ast}_{i}(t)) - \sum_{j\in\mathcal{B}(i)}x^{\ast}_{j}(\sigma^{\ast}_{j}(t)) 
	\right)
	\mathrm{d}t &=
	\int_{\ClT} 
	\left(
	x_{i}(\sigma^{\ast}_{i}(t)) \dot{\sigma}^{\ast}_{i}(t)- \sum_{j\in\mathcal{B}(i)}x^{\ast}_{j}(\sigma^{\ast}_{j}(t)) \dot{\sigma}^{\ast}_{j}(t)
	\right)
	\mathrm{d}t 
	\\
	&=
	\left[ 
		D_{i}(\sigma^{\ast}_{i}(t)) - \sum_{j\in\mathcal{B}(i)}D_{j}(\sigma^{\ast}_{j}(t))
	\right]_{\ClT}
	\\
	&= Q_{i}.
\end{align}

We finally confirm that Condition~\eqref{Eq_DSO-MPCC_cnst_q_Nonneg} holds.
To this end, we look at the mathematical expression of the inflow $q^{*}_{i}(t)$ (Eq.~\eqref{Eq:Trans_q}) in detail by dividing the situation into the following two cases: (i) $A_{i}(\sigma^{\ast}_{i}(t)) \neq D_{i}(\sigma^{\ast}_{i}(t))$ and (ii) $A_{i}(\sigma^{\ast}_{i}(t)) = D_{i}(\sigma^{\ast}_{i}(t))$.

Consider first Case (i).
This condition means that in an original traffic state, there exists a queue on link $i$ at the time $\sigma_{i}^{*}(t)$, and thus $x_{i}(\sigma_{i}^{*}(t)) = \mu_{i}$.
Substituting this into Eq.~\eqref{Eq:Trans_q}, we have 
\begin{align}
    q^{\ast}_{i}(t) = \mu_{i} - \sum_{j\in\mathcal{B}(i)}x_{j}(\sigma^{\ast}_{j}(t)).
\end{align}
Then, from Eq.~\eqref{Eq_Asm_Capa}, we see that $q^{\ast}_{i}(t) \geq 0$ in this case.

Consider next Case (ii).
This condition means that a queue does not exist at that time.
In this case, we first derive the following equation by combining the condition and the flow conservation condition:
\begin{align}
    &
    D_{i}(\sigma^{\ast}_{i}(t)) = A_{i}(\sigma^{\ast}_{i}(t)) 
    = Q_{i}(\sigma^{\ast}_{i}(t)) + 	\sum_{j\in\mathcal{B}(i)} D_{j}(\sigma^{\ast}_{j}(t)),
    &&
    \forall i \in \ClN, \quad \forall t \in \ClT,
		\label{Eq_q_Nonneg_cumulative}
\end{align}
where $Q_{i}(\sigma^{\ast}_{i}(t))$ represents the cumulative inflow of $(i)$-commuters by time $\sigma^{\ast}_{i}(t)$.
Differentiating Eq.~\eqref{Eq_q_Nonneg_cumulative} with respect to $t$ yields 
\begin{align}
	&\diff{Q_{i}(\sigma^{\ast}_{i}(t))}{\sigma^{\ast}_{i}(t)} \dot{\sigma}^{\ast}_{i}(t)
	+
	\sum_{j\in\mathcal{B}(i)} x_{j}(\sigma^{\ast}_{j}(t)) \dot{\sigma}^{\ast}_{j}(t)= x_{i}(\sigma^{\ast}_{i}(t))\dot{\sigma}^{\ast}_{i}(t)&&\\
	&\Leftrightarrow\quad 
	\diff{Q_{i}(\sigma^{\ast}_{i}(t))}{\sigma^{\ast}_{i}(t)}
	+
	\sum_{j\in\mathcal{B}(i)} x_{j}(\sigma^{\ast}_{j}(t))= x_{i}(\sigma^{\ast}_{i}(t))
	&&(\because \dot{\sigma}^{\ast}_{i}(t) = \dot{\sigma}^{\ast}_{j}(t) = 1).
\end{align}
We thus have the following relationship:
\begin{align}
    q^{\ast}_{i}(t) = 
    x_{i}(\sigma^{\ast}_{i}(t)) - \sum_{j\in\mathcal{B}(i)} x_{j}(\sigma^{\ast}_{j}(t)) 
	= \diff{Q_{i}(\sigma^{\ast}_{i}(t))}{\sigma^{\ast}_{i}(t)} \geq 0.
\end{align}
Hence, Condition~\eqref{Eq_DSO-MPCC_cnst_q_Nonneg} holds.
This completes the proof of the lemma.

We are now ready to formally prove the non-existence theorem based on this lemma, by showing that the new traffic state resulting from the lemma always decreases the total system cost.
From Eq.~\eqref{Eq:Trans_q}, we derive the following equation:
\begin{align}
    &\sum_{i \in \ClN} q^{\ast}_{i}(t) = \sum_{j\in\mathcal{B}(0)}x_{j}^{*}(\sigma_{j}^{*}(t))
    = \sum_{j\in\mathcal{B}(0)}x_{j}(\sigma_{j}^{*}(t))
	&&\forall t \in \ClT.
\end{align}
By definition, for a node $j\in\mathcal{B}(0)$ (just upstream of the unique destination), $\sigma_{j}^{*}(t) = \sigma_{j}(t) = t$.
This means that: 
\begin{align}
    &\sum_{i \in \ClN} q^{\ast}_{i}(t) = \sum_{j\in\mathcal{B}(0)}x_{j}(\sigma_{j}(t))
    = \sum_{i \in \ClN} q_{i}(t),
    &&\forall t \in \ClT.
\end{align}
Therefore, we have the following relationship between the objective functions $f(\Vtq, \Vtw)$ and $f(\Vtq^{\ast}, \Vtw^{\ast})$:
\begin{align}
    f(\Vtq^{\ast}, \Vtw^{\ast}) - 
	f(\Vtq, \Vtw)
	&= 
	\sum_{i \in \ClN}
	\int_{\ClT} 
		c(t) 
	q^{\ast}_{i}(t) \mathrm{d} t
	-
	\sum_{i \in \ClN}
    \int_{\ClT} 
    \left(
    c(t) + \sum_{j;j \in \underline{\ClN}(i)} w_{j}(t)
    \right)
    q_{i}(t) \mathrm{d} t
	\\
	&\int_{\mathcal{T}}c(t)\left(\sum_{i\in\mathcal{N}}q^{\ast}_{i}(t) - 
	    \sum_{i\in\mathcal{N}}q_{i}(t)
	\right)\mathrm{d}t
	- \sum_{i \in \ClN} \int_{\ClT} w_{i}(t)\sum_{j;j \in \overline{\ClN}(i)} q_{j}(t) \mathrm{d} t\\
	&=
	- \sum_{i \in \ClN} \int_{\ClT} w_{i}(t)\sum_{j;j \in \overline{\ClN}(i)} q_{j}(t) \mathrm{d} t \leq 0.
\end{align}
\noindent Note that the equality holds if and only if $\Vtw = \Vt0$.

Therefore, transforming a traffic state with queues by \textbf{Lemma~\ref{Lemm:DeadWeight}}, we can derive the traffic state without queues where the total system cost is always lower than the cost in the original traffic state.
This completes the proof of the theorem.\qed

\section{Concluding remarks}
This study explored the non-existence property in the DSO state in a directed rooted tree network with a many-to-one origin destination pattern.
We first straightforwardly formulated a DSO problem as an MPCC that minimized the schedule and queueing delay costs with queueing conditions.
We then showed that for an arbitrary traffic state with queues, we can construct another traffic state where no queues exist without changing the total schedule delay costs.
This proves that the queueing delay costs are deadweight loss and thus they must not exist in a DSO state.
Additionally, as an application of the non-existence property, we showed that the DSO problem can be formulated as an LP.
In other words, we found there exists the possibility of a transformation from an MPCC to an LP under a certain condition.

As part of our future work, we will investigate whether the non-existence theorem of queues holds in more general settings, such as in networks with multiple routes and heterogenous commuters.
If the non-existence theorem does not hold generally, it would be important to clarify the conditions under which the theorem does hold.
Moreover, our results would be useful for analyzing the relationship between a DSO state without queues and a dynamic user equilibrium (DUE) state, in which each user optimizes his/her own utility and thus queues inevitably occur.
In fact, recent studies~\citep{Fu2021-qg,SAKAI2021-th,Sakai2022-ab} showed how the flow and cost patterns in the DUE state correspond to those in the solution to the [DSO-LP] problem.
Therefore, it would be interesting to investigate the relationship between the DSO and DUE states from the perspective of the non-existence theorem.


\section*{Acknowledgments}
This work was supported by JSPS KAKENHI Grant Numbers JP20J21744, JP20H02267, and JP21H01448.

\bibliography{myrefs}

\end{document}